\documentclass[12pt]{amsart}

\usepackage{amsmath,amsfonts,amsthm,amsopn,cite,mathrsfs}
\newtheorem{tm}{Theorem}
\newtheorem{lemma}[tm]{Lemma}
\newtheorem{prop}[tm]{Proposition}
\newtheorem{cor}[tm]{Corollary}
\theoremstyle{definition}
\newtheorem{definition}{Definition}

\newcommand{\field}[1]{\mathbb{#1}}
\newcommand{\bR}{\field{R}}        
\newcommand{\bN}{\field{N}}        
\newcommand{\bZ}{\field{Z}}        
\newcommand{\bC}{\field{C}}        
\def\rd{\bR^d}

\def\inv{^{-1}}

 \def\cH{\mathcal{H}}
 \def\cB{\mathcal{B}}

 \def\cM{\mathcal{M}}

 \def\cA{\mathcal{A}}

 \def\cC{\mathcal{C}}
 
 \def\cO{\mathcal{O}}

 \def\zd{\bZ^d}
 \def\lpm{\ell ^p_m}
\newcommand{\fif}{if and only if}

\newcommand{\ant}{\widetilde{A_n}}

\newcommand{\lz}{\ell ^2 (\bZ ^d )}
\newcommand{\avo}{\cA ^1_v}
\newcommand{\av}{\cA _v}
\newcommand{\as}{\cA _s}

\newcommand{\bn}{\mathbb{N}}

\newcommand{\bz}{\mathbb{Z}^d}

\newcommand{\Cal}{\mathcal}
\newcommand{\lan}{\langle}
\newcommand{\ran}{\rangle}
 \newcommand{\TCr}{\widetilde{D_{r,n}}}
 \newcommand{\TEr}{\widetilde{E_{r,n}}}
 \newcommand{\TBn}{\widetilde{B_{n}}}
\newcommand{\TCri}{\widetilde{D_{r,n}}^{-1}}

\newcommand{\Cri}{C_{r,n}^{-1}}
\newcommand{\Dri}{D_{r,n}^{-1}}

\newcommand{\el}{{\ell}^2(\bz)}
\newcommand{\ac}{\Cal A}

\def\Cst{\mathbb C}
\def\Zst{\mathbb Z}
\def\Rst{\mathbb R}
\def\AT{A^{\ast}}
\def\at{a^{\ast}}

\def\Arn{A_{r,n}}
\def\xrn{x_{r,n}}

\def\tI{\operatorname{I}}
\def\tII{\operatorname{II}}
\def\tIII{\operatorname{III}}

\def\phi{\varphi}

\theoremstyle{remark}
\newtheorem{remark}{Remark}[section]
\numberwithin{equation}{section}

\begin{document}

\title{Quantitative Estimates for the Finite Section Method}
\author{Karlheinz Gr\"ochenig}
\author{Ziemowit Rzeszotnik}
\address{Faculty of Mathematics \\
University of Vienna \\
Nordbergstrasse 15 \\
A-1090 Vienna, Austria}
\email{karlheinz.groechenig@univie.ac.at,ziemowit.rzeszotnik@univie.ac.at}

\author{Thomas Strohmer}
\address{Department of Mathematics, University of California, Davis,
  CA 95616-8633, USA}
\email{strohmer@math.ucdavis.edu}
\date{}
\thanks{K.~G.~and Z.~R. were  supported by the Marie-Curie Excellence
  Grant MEXT-CT 2004-517154, T.~S.\ was partially supported by NSF DMS
Grant 0511461.} 
\maketitle

\begin{abstract}
The finite section method is a classical scheme to approximate the solution
of an infinite system of linear equations. We present quantitative
estimates for the rate of the convergence of the finite section method
on weighted $\ell ^p$-spaces.
Our approach uses recent results from the theory of Banach algebras
of matrices with off-diagonal decay. Furthermore, we demonstrate that
Banach algebra theory provides a natural framework for deriving a finite
section method that is applicable to large classes of non-hermitian matrices.
An example from digital communication illustrates the practical
usefulness of the proposed theoretical framework.
\end{abstract}

\section{Introduction}

Many of the concrete applications of mathematics in science and
engineering eventually result in a problem involving linear operator equations.
This problem can be usually represented as a linear system of 
equations (for instance by discretizing an integral equation or because 
the operator equation is already given on some sequence space) of
the form
\begin{equation}
\label{Axb}
Ax  = b,
\end{equation}
where $A$ is an infinite matrix $A = (a_{kl})_{k,l\in \bZ }$ and 
$b$ belongs to some  Banach space of sequences. 
Solving  linear equations with infinitely many variables is a problem
of functional analysis, while solving  equations with finitely many variables
is one of the main themes of linear algebra. Numerical analysis bridges
the gap between these  areas.
A fundamental problem of numerical analysis is thus  to find  a
finite-dimensional model for~\eqref{Axb} whose solution approximates the solution of the 
original infinite-dimensional  problem with any desired accuracy. This problem often leads
to delicate questions of stability and convergence. 

A simple and useful approach is the {\em finite-section 
method}~\cite{GF74,HRS01}. Let 
$$P_n b = ( \dots ,  0,  b_{-n}, b_{-n+1}, \dots ,
b_{n-1}, b_n, 0, \dots )$$ 
be the orthogonal projection onto a $2n+1$-dimensional subspace
. We set
\begin{equation}
\label{Axb1}
A_n = P_n A P_n \quad \quad \text{ and } \quad \quad b_n = P_n b \, ,
\end{equation}
and try to solve the finite system 
\begin{equation}
\label{Axb2}
A_n x_n = b_n  
\end{equation}
for properly chosen $n$. The crucial question is then:
What is the relation between the numerical solution $x_n$ and
the actual solution $x$?

This problem has been analyzed in depth for the case
of convolution operators and Toeplitz matrices in the pioneering
work of Gohberg, e.g. see~\cite{GF74}. Important generalizations and 
extensions in the Toeplitz setting can be found in~\cite{BS83,BS90}.
Rabinovich et al.\ derive necessary and sufficient conditions for the 
convergence of the finite section method in terms of the so-called 
limit operator~\cite{RRS04}, which does not necessarily require any
Toeplitz structure. These conditions, while intriguing, are not
always easy to verify in practice.

A general theory for the approximation by  finite-section is based on  the
powerful methods  of $C^{\ast}$-algebras and  has been developed by B\"ottcher,
Silbermann, and coworkers, see for instance~\cite{BS83,HRS01}.
Their framework leads to many attractive and deep results about 
the applicability of the finite section method as well as other approximation 
methods. William Arveson goes a step further and concludes
 that  {\em ``numerical problems involving infinite dimensional operators require 
a reformulation in terms of $C^{\ast}$-algebras''}~\cite{Arv94}. However, $C^{\ast}$-algebras 
have some limitations. It  was already pointed out in~\cite{HRS01}
that $C^{\ast}$-algebra techniques  do not  yield any information about the speed of
convergence of the finite section method. 
An answer to this question 
is obviously not only of theoretical interest,  but it is   important
for real applications. For instance, we want to choose $n$ in~\eqref{Axb2} large enough
to get a sufficiently accurate solution, but on the other hand, $n$
should be small enough to bound the computational complexity which in
general is of order $\cO (n^3)$. 
Theorems about the speed of convergence will give a quantitative indication
for how increasing $n$ will impact the accuracy of the solution.
Some results about the speed of convergence for the special case of Toeplitz 
matrices can be found in~\cite{Str98a,Str00,RRS01,GGK03}. 
In~\cite{GGK03} the 
convergence in the $\ell^p$-norm ($1\le p<\infty$) is analyzed.

In this paper we present a thorough analysis  of the  convergence of 
the finite section method for positive definite matrices as well as 
for non-hermitian ones. Specifically, we solve the following
problems. 

(a) We study the finite section method on weighted $\ell
^p$-spaces. If the input vector $b$ belongs to a weighted space $\ell
^p_m$, then, under suitable assumptions on the matrix 
$A$,  the finite section method converges in the norm of $\lpm $. 

(b) We obtain quantitative estimates for the rate of convergence of
$x_n$ to $x$ in various weighted $\ell^p$-norms.  

(c) We define a modified version of finite sections, the
\emph{non-symmetric finite section method},  and show that
this  method converges also for non-symmetric matrices. The finite
section method for non-symmetric matrices raises a number of
rather difficult questions and has motivated a large part of
\cite{HRS01}. Even for the classical case of Laurent operators
(Toeplitz matrices) our  approach enlarges considerably the class of
matrices to which the finite section method can be applied.  

As we work with Banach spaces of sequences, the methods will be taken
from the theory of $B^*$-algebras (involutive Banach algebras) instead
of $C^*$-algebras which suit only  Hilbert spaces.  
The key property of the matrices $A$ is their off-diagonal decay; we
will rely heavily on  recent results from 
the theory of Banach algebras of matrices. In fact, an important
technical part of our analysis is to establish a finite section
property of infinite-dimensional matrix algebras. 

The paper is organized as follows. In Section~2 we recall the well
known proof for the convergence of the finite section method for
positive invertible matrices and take it as a model for more general
statements. In Section~3 we introduce several Banach algebras of
infinite matrices and collect their fundamental properties. Section~4
is devoted to the notion of inverse-closedness and spectral invariance
in Banach algebras and their relation to the finite section method. In
Section~5 we establish the convergence of the finite section method on
weighted $\ell ^p$-spaces, in Section~6 we derive quantitative
estimates. In Section~7 we investigate a version of the finite section method for
non-symmetric matrices, and in the final Section~8 we briefly discuss
an application to wireless communications.

\section{Convergence of the finite section method}

It is well known that for positive definite matrices the 
finite section method works in principle, see, e.g.,~\cite{HRS01}. The proof is 
instructive and exhibits what is necessary for an understanding of the 
finite section method.

Recall that if $\cA $ is an algebra, then  the spectrum of an element
$A \in \cA $ is defined to be the set $\sigma _{\cA } (A) = \{ \lambda
\in \bC : (A-\lambda I) \, \mathrm{is \, not \,  invertible }\}$. If the
algebra is $\cB (\cH )$, the bounded operators on some Hilbert space,
we usually omit the reference to the algebra and write simply $\sigma
(A)$ for the spectrum. For self-adjoint operators on $\cH $ we denote
the extremal spectral values of $\sigma (A)$ by $\lambda_-=\min\sigma(A)$ and
$\lambda_+=\max\sigma(A)$, so that $\sigma (A) \subseteq [\lambda _-,
\lambda _+]$. 

We will analyze the finite section method for multidimensional index sets
of the form $\Zst^d$. 
To that end we define the projection $P_n$ in dimension $d>1$.
We set $C_n = [-n,n]^d \cap \zd $,  the integer vectors in the cube of
length $2n$ centered at the origin. Then the projection $P_n$ is defined by 
$(P_ny)(k) = \chi _{[-n,n]^d}(k) y(k) = \chi _{C_n}(k) y(k)$ for $k\in
\zd $. 
The range of $P_n$ is a subspace of $\ell ^2 (\zd ) $ of dimension $(2n+1)^d$
and will be identified with $\bC ^{(2n+1)^d}$. The finite section is
then defined to be $A_n=P_nAP_n$. 
By definition, $A_n$ is  a (finite rank) operator acting on
$\ell ^2 (\zd ) $, but we often interpret $A_n$ as a finite $(2n+1)^d \times (2n+1)^d$-matrix
acting on $\bC ^{(2n+1)^d}$. In particular,  by $A_n \inv $ we
understand the inverse of this finite matrix, but clearly 
$A_n $ cannot be invertible  on $\ell ^2 (\zd ) $.

We mention  that our results could also be formulated
with respect to other index sets. 
\begin{tm} \label{cstar}
  If $A$ is a positive and (boundedly) invertible operator on $\lz $,
  then $x_n $ converges to $x$ in $\lz $. 
\end{tm}

\begin{proof}
\textbf{Step 1.} Since by hypothesis, $\sigma (A) \subseteq [\lambda_-,
\lambda _+]\subseteq (0,\infty )$, we have 
$$
\lambda _- \|P_n b\|_2^2 \leq \langle AP_n b , P_n b \rangle = \langle
A_n b, b \rangle \leq \lambda _+ \|P_n b \|_2^2 \, .
$$
Consequently  on the invariant subspace $P_n \lz \simeq \bC ^{(2n+1)^d}$ 
$$
\sigma (A_n ) \subseteq [\lambda
  _-, \lambda _+]
$$
independent of $n$. In particular, each $A_n $ is invertible on $\bC
^{(2n+1)^d}$ and
\begin{equation}
  \label{eq:fe1}
  \sup _{n\in \bN } \|A_n \inv \|_{op} \leq \lambda _-\inv = \|A\inv
  \|_{op} \, .
\end{equation}

\textbf{Step 2.} Define an extension of  $A_n $ by 
\begin{equation}
  \label{eq:1}
  \ant = A_n + \lambda _+ (I-P_n ) \, .
\end{equation}
 Then 
$\sigma (\ant ) \subseteq [\lambda _-, \lambda _+ ]$, and  
 all matrices $\ant $ are
invertible on $\lz $. Furthermore, $\ant \inv = A_n \inv + \lambda _+ \inv
(I-P_n )$ and $\ant $ converges to $A$ in the strong operator
topology.

\textbf{Step 3.} (Lemma of Kantorovich). Since
\begin{eqnarray}
\|\ant \inv b - A\inv b \|_2 &=&\| \ant \inv (A-\ant ) A\inv b \|_2
\notag  \\
&\leq & \sup _n \|\ant \inv \|_{op} \, \| (A-\ant ) A\inv b \|_2  \, ,
\end{eqnarray}
the strong convergence $\ant \rightharpoonup A $ implies that $\ant
\inv $ converges strongly to $A\inv $. 

\textbf{Step 4.} Recall $A_n x_n = b_n $ and $Ax=b$. Then
\begin{eqnarray}
  \|x-x_n \|_2 &=& \|A\inv b - A_n \inv b_n \|_2 = \| A\inv b - A_n
  \inv P_n  b \|_2 \notag  \\
&\leq &\| (A\inv - \ant \inv ) b\|_2 + \| \ant \inv (b - P_n b) \|_2 = 
\tI + \tII \, .
\end{eqnarray}
The first term goes to zero by Step 3, and the second term
is estimated by 
$$ \tII \leq \sup _n \|\ant \inv \|_{op} \, \|b-P_n b\|_2
\leq \lambda _- \inv \|b-P_n b\|_2$$
and also goes to zero. 
\end{proof}

The above theorem uses the $\lz $-norm, so this is the realm of 
$C^*$-algebra techniques, cf.\ the work of B\"ottcher, Silbermann, 
et al.~\cite{BS83,HRS01}.

Several questions arise naturally in the context of the finite section
method:\\
1. Does the finite section method also converge in other norms, e.g.,
in weighted $\ell ^p$-norms? \\
2. Can we derive quantitative estimates? If the finite section method works,  how fast
does $x_n $ converge to $x$? What conditions on the matrix $A$ and
the input vector $b$  are
required to quantify the rate of convergence $x_n \to x$?\\
3. What conditions and modifications are required (if any) to make the
finite section method work for matrices that are not hermitian? 

For an  answer of the first question, we make the following
observation:  The simple argument above extends almost word by word,
provided we can show  the 
following properties: 
\begin{itemize}
\item[(1)]
 Both $A$ and $A\inv $ are bounded on $\lpm $,  
\item[(2)] 
 $\sup _n \|\ant \inv \|_{\ell ^p_m \to \ell ^p_m}$ is
finite, and
\item[(3)]
 the finite sequences are dense in $\lpm $. 
\end{itemize}
The answers to the other two questions also revolve around the above
observation as well as on properties of certain involutive Banach algebras,
which will be introduced in the next section.

\section{A Class of Banach Algebras of Matrices}

To understand the asymptotic behavior of the finite section method on 
Banach spaces,
we need to resort to  Banach algebra methods. We first  consider some
typical matrix norms that express various forms of off-diagonal decay.  
Our approach is partly motivated by some forms of  off-diagonal decay
that is  observed in various applications, such as signal and image 
processing, digital communication, and quantum physics. A different way
of describing off-diagonal decay of matrices (and operators) is given
by the notion of band-dominated operators~\cite{RRS01}.

\textbf{Weights.} Off-diagonal decay is quantified  by means of weight
functions. 
 A non-negative function $v$ on $\zd $ is called an {\em admissible weight}
if it satisfies the following properties:
\begin{itemize}
\item[(i)] $v$ is even  and normalized such
that $v(0) = 1$.
\item[(ii)] $v$ is submultiplicative, i.e., $v(k+l)\le v(k) v(l)$ for
  all $k,l\in \zd$.
\end{itemize}
The assumption that $v$ is even assures that the corresponding Banach
algebra is closed under taking the adjoint $A^*$. The weight $v$ is said to 
satisfy the \emph{Gelfand-Raikov-Shilov (GRS) condition}~\cite{GRS64}, if  
\begin{equation}
\underset{n \to \infty}{\lim} v(n k)^{\frac{1}{n}} = 1 \qquad
\text{for all $k\in \zd  $}.
\label{GRS}
\end{equation}
This property is crucial for the inverse-closedness of Banach
algebras, see Theorem~\ref{inverseclosed} below.
The standard weight functions on $\zd $ are of the form
$$v(x) = e^{a d (x)^b} (1+d(x))^s \, , 
$$ 
where $d(x)$ is a norm on $\rd $. Such a weight is submultiplicative, when  
$a,s \geq 0$ and $0\leq b \leq 1$; $v$  satisfies the GRS-condition,
\fif\ $0\leq b < 1$.

Consider the following conditions on matrices.

1. \emph{The Jaffard class} is  defined by polynomial decay off the
diagonal. Let $\cA _s$ be the class of matrices $A= (a_{kl}), k,l \in
\zd$, such that
\begin{equation}
  \label{eq:m5}
  |a_{kl}| \leq C (1+|k-l|)^{-s} \quad \quad \forall k,l \in \zd 
\end{equation}
with norm $\|A\|_{\cA _s} = \sup _{k,l \in \zd } |a_{kl} |
(1+|k-l|)^{s}$. 


2. \emph{More general off-diagonal decay.}
Let $v$ be an admissible   weight on $\zd $ that satisfies the
following additional conditions: $v\inv \in \ell ^1(\zd )$ and $v\inv
\ast v\inv \leq C v\inv $ ($v$ is called \emph{subconvolutive}).  We define the
  Banach space  $\cA _v $ by the  norm
  \begin{equation}
    \label{eq:s12}
    \|A\|_{\cA _v } = \sup _{k,l \in \zd}     |a_{kl}| v(k-l) \, ,
  \end{equation}

3. \emph{Schur-type conditions.}
Let $v$ be an admissible weight. The class  $\avo $ consists of   all
matrices $A = (a_{kl})_{k,l \in 
  \zd}$ such that
  \begin{equation}
    \label{eq5}
        \sup_{k\in \zd} \sum _{l\in \zd} |a_{kl}| \, v(k-l) < \infty
\quad \text{ and } \quad      \sup_{l\in \zd} \sum _{k\in \zd} |a_{kl}| \,
v(k-l) < \infty
  \end{equation}
with norm
\begin{equation}
  \label{eq7}
  \| A \|_{\avo } = \max \big\{ \sup _{k\in \zd } \sum _{l \in \zd} |a_{kl}|
  v(k-l) \, , \, \sup _{l\in \zd } \sum _{k \in \zd} |a_{kl}|
  v(k-l)\big\} \, .
\end{equation}

4. \emph{The Gohberg-Baskakov-Sj\"ostrand class.} For any admissible
weight $v$ we define the class  $\cC _v $
as the space of
  all matrices $A = (a_{kl})_{k,l \in \zd }$ such that the norm
  \begin{equation}
    \label{eq:17}
    \|A\|_{\cC _v} := \sum _{l\in \zd } \sup _{k\in \zd } |a_{k,k-l}|\, v(l)
       \end{equation}
is finite.  An alternative way to define the norm on $\cC _v$ is
\begin{equation}
\|A\|_{\cC _v} = \inf \{ \|\alpha \|_{ \ell ^1 _v } : |a_{kl} | \leq
\alpha (k-l) \} \, .
\label{eq:17a}
\end{equation}

5. A further generalization is due to Sun~\cite{Sun05}. Roughly speaking, Sun's
class amounts to an interpolation between $\cC _v $ and $\cA _v$ or
between $\avo $ and $\av $. Our results also hold for Sun's class, but
to avoid a  jungle of indices, we stick to the simple classes
defined above and leave the reformulation of our results in Sun's case
to the reader.

These Banach spaces of matrices have the following elementary properties.

\begin{lemma} \label{bound}
  Let $v$ be an admissible weight and $\cA $ be one of the algebras
  $\cA _s$ for $s>d$, $\cA _v , \avo , \cC _v 
  $.  Then $\cA $ has the following   properties: 

(a) Both $\cA _v^1$ and $\cC _v$ are involutive Banach algebras (i.e.,
$B^{\ast}$-algebras) with the norms defined in~\eqref{eq5} and \eqref{eq7}. 
$\cA _v$ and $\cA _s, s>d$ can be equipped with an equivalent norm so 
that they become involutive Banach algebras. 

(b) If $A \in \cA $, then $A $ is bounded on $\lz $. 

(c) If $A\in \cA $ and  $|b_{kl}|\leq |a_{kl}|$ for all $k,l \in \zd
$, then $B\in \cA $ and $\|B\|_{\cA } \leq \|\cA \|_{\cA}$. ($\cA $ is
a \emph{solid} algebra). 
\end{lemma}

\begin{proof}
Properties (a) and (c) are easy and follow directly from the
definition of the matrix norms. The statements about $\cA_s$ 
and $\cA_v$ are proven in~\cite{GL04}.
(b) is a consequence of  Schur's test.
\end{proof}

Next we study the spectrum of matrices belonging to one of these
Banach algebras. 
\begin{definition}\label{definverse}
We say that $\ac$ is inverse-closed in $\cB (\ell ^2(\zd ))$, 
if for every $A\in \ac$ that is invertible on $\el$ we have that $A^{-1}\in \ac$.
\end{definition} 

Our next theorem states that the matrix algebras introduced above are
inverse-closed as long as $v$ satisfies the  GRS-condition. The precise formulation is slightly  more complicated,
because we need to be a bit pedantic about  the weights. 

\begin{tm}[Inverse-closedness]
\label{inverseclosed}
 Let $v$ be an admissible weight that  satisfies the
 GRS-condition, i.e., $\lim _{n\to \infty }   v(nk)^{1/n} = 1$ for all $k\in \zd $.

(a) Assume that $v\inv \in \ell ^1(\zd )$ and $v\inv
\ast v\inv \leq C v\inv$,
then   $\av $ is  inverse-closed in $\cB (\lz )$. 
In particular $\cA _s$ for $s>d$ possesses this property. 

(b) If $v(k) \geq C (1+|k|)^\delta $ for some $\delta >0$, then  $\avo
$ is inverse-closed in  $\cB (\lz )$. 
 
(c) $\cC _v$ is inverse-closed in  $\cB (\lz )$ for arbitrary
admissible weights with the GRS-property. 
\end{tm}

\begin{remark} 
The inverse-closedness is the key property and lies rather
deep. While for $C^{\ast}$-(sub)algebras this property is inherent, for Banach
algebras it is always hard to prove. Inverse-closedness for $\cA _s$ is
due to Jaffard~\cite{Jaf90} and Baskakov~\cite{Bas90,Bas97}, a simple proof 
is given in ~\cite{Sun05}. For $\av $ it was
proved by Baskakov~\cite{Bas97} and reproved in a different way in
~\cite{GL04}. The result for $\cC _v$ with $v\equiv 1$  is due to
Gohberg-Kasshoek-Wordeman~\cite{GKW89} and was rediscovered by
Sj\"ostrand~\cite{Sjo95}, the case of arbitrary weights is due to
Baskakov~\cite{Bas97}, the algebra $\avo $ was treated by one of us with
Leinert~\cite{GL03}. More general conditions were announced by
Sun~\cite{Sun05}. 
\end{remark}

The following properties are well-known consequences of
inverse-closedness. 

\begin{cor}[Spectral invariance] \label{closedgraph}
Let $\cA $ be one of the algebras $\cA _s$, $\av $, $\avo $, or $\cC
_v$ and assume that $v$ satisfies the conditions of
Theorem~\ref{inverseclosed}. 
Then \\
(a)  $\sigma _{\cA } (A) = \sigma (A) $ (the spectrum in the algebra
$\cA $ coincides with the spectrum of $A$ as an operator on $\lz$) 
  
(b) If $A$ is bounded on $\lpm $ for all $A \in \cA $, then the
operator norm satisfies
\begin{equation}
\label{lpbound}
\|A\|_{\ell ^p_m \to \ell ^p_m} \leq  C \| A
\|_{\cA } \quad \text{for all $A \in \cA$,}
\end{equation}
and 
$$
\sigma _{\ell ^p_m } (A) \subseteq \sigma (A) 
$$
(the spectrum is almost independent of the
space $A$ acts on).
\end{cor}
\begin{remark}
Statement (a) is equivalent to inverse-closedness, the norm
estimate in (b) follows from the closed graph theorem, the inclusion
of the spectra is an immediate consequence of
Theorem~\ref{inverseclosed}. 
\end{remark}

Let us emphasize that in our analysis of the finite section method we
only need that the algebra $\cA $ acts boundedly on $\lpm $. In order
to understand how the weight $m$ depends on the submultiplicative weight
used to parametrize the off-diagonal decay, let us briefly discuss some
sufficient conditions for the bounded action of $\cA $ on $\lpm $. The
weights $m$ satisfy slightly different conditions.
Let $v$ be an admissible weight. The class of $v$-moderate weights is  
\begin{equation}
\cM _v = \Big\{m\ge 0: \underset{k\in \zd}{\sup} 
\frac{m(k+l)}{m(k)} \le Cv(l), \quad \forall \,l  \in \zd \Big\}.
\label{moderateweights}
\end{equation}
For example, if $a, s \in \bR$ are  arbitrary, then $m(x) = e^{a d (x)^b} (1+d(x))^s$ is $e^{|a|
 d (x)^b} (1+d (x))^{|s|}$-moderate. 

The explicit examples of Banach algebras  discussed above 
all act on the entire range of $\lpm $ for
$1\leq p \leq \infty $ and a family of moderate weights associated to
$v$. The following lemma provides some explicit sufficient  conditions 
on $m$ for $\cA _v$, $\cA ^1_v$ or $\cC _v$ to act boundedly on $\lpm $.

\begin{lemma}\label{lembound}
Let $v$ be an admissible weight. 

 (a) If $A \in \avo $, then $A$ is bounded simultaneously on all $\ell
 ^p_m (\zd )$ for $1 \leq  p \leq \infty $ and $m \in \cM _v$. 

(b) If $A \in \av$ and $v_0(k) :=v(k)/ (1+|k|)^s $ is
submultiplicative for some $s>d$, then  $A$ is  bounded simultaneously on all $\ell
 ^p_m (\zd )$ for $1 \leq p \le \infty $ and $m \in \cM _{v_0}$.

(c) If  $A \in \cA _v$, then $A$ is bounded on $\ell ^\infty _v (\zd )$.   

(d) If $A\in \cC _v$, then $A$ is bounded on all $\ell ^p _m (\zd )$
for $1 \leq  p \leq \infty $ and $m \in \cM _v$. 
\end{lemma}

\begin{proof}
For completeness we sketch the easy proof.

(a) First, let $p=1$, $c\in \ell^1_m(\zd)$ and $A\in\avo$. Then, since $m(k)\le C v(k-l)m(l)$, we obtain
\begin{align*}
\|Ac\|_{\ell^1_m} & =  \sum_{k\in\zd} \Big|\sum_{l\in\zd} a_{kl} c_l \Big|m(k) \le
 C \sum_{k\in\zd} \sum_{l\in\zd} |a_{kl}|\, |c_l| v(k-l)m(l)  \\ 
& \le C\sum_{l\in\zd} \Big(\sup_{l\in\zd}\sum_{k\in\zd} |a_{kl}| v(k-l)\Big)|c_l|m(l)=C\|A\|_{\avo}
\|c\|_{\ell^{1}_{m}}.
\end{align*}
Next, let $p=\infty$ and $c\in\ell^{\infty}_{m}$. Then, as before
\begin{align*}
\|Ac\|_{\ell^\infty_m} & =  \sup_{k\in\zd} \Big|\sum_{l\in\zd} a_{kl} c_l \Big|m(k) \le
 C \sup_{k\in\zd} \sum_{l\in\zd} |a_{kl}|\, |c_l| v(k-l)m(l)  \\ 
& \le C\Big(\sup_{l\in\zd} |c_l|m(l)\Big)\sup_{k\in\zd}\sum_{l\in\zd} |a_{kl}| v(k-l)=
C\|A\|_{\avo}\|c\|_{\ell^{\infty}_{m}}.
\end{align*}
The boundedness on $\ell^p_m(\zd)$ for $1<p<\infty$ now follows by interpolation.

(b) and (d) follow from the easy embeddings $\av\hookrightarrow\cA_{v_0}^1$, 
$C_v\subseteq \avo$ and from (a).

(c) uses the subconvolutivity of $v$. Let $A\in\av$ and
$c\in\ell^\infty_v(\zd)$. Then, $|a_{kl}|\le \|A\|_{\av}v(k-l)^{-1}$
and $|c_l|\le \|c\|_{\ell_v^\infty}v(l)^{-1}$. Consequently, 
\begin{align*}
\|Ac\|_{\ell^\infty_v} & =  \sup_{k\in\zd} \Big|\sum_{l\in\zd} a_{kl} c_l \Big|v(k) \\ &\le \|A\|_{\av}\|c\|_{\ell^{\infty}_{v}}
\sup_{k\in\zd} \sum_{l\in\zd} \frac1{v(k-l)}\frac1{v(l)}v(k)\le C\|A\|_{\av}\|c\|_{\ell^{\infty}_{v}},
\end{align*}
because $(v^{-1} \ast v^{-1})(k)\le C{v(k)^{-1}}$.
\end{proof}



The matrices of the Banach algebras introduced above can be 
considered as approximate banded matrices. This might suggest that it
would be sufficient to set those entries smaller than some threshold to zero 
and simply work with banded matrices, which are a special case of sparse
matrices. At first sight this may seem appealing, since invertible 
banded matrices have inverses with exponentially fast off-diagonal 
decay~\cite{DMS84}. However there is an important difference. 
In many applications, cf.~\cite{Jaf90,Str98a,Str00,Str05} 
the matrix entries do decay off the diagonal, but thresholding would 
still leave us with banded matrices with the number of
non-zero diagonals easily in the order of several dozens.
The theoretical prediction for the  decay of the inverse of such banded 
matrices is so slow that it is meaningless for practical purposes. The reason 
is that by resorting to banded matrices we have neglected the decay of the 
entries above the chosen threshold. Thus banded matrices
are simply not the most suitable model to capture the decay behavior of those
matrices and their inverses.

\section{Finite Sections in Matrix Algebras}\label{zioma}

We first  study the finite sections of matrices belonging to
an inverse-closed matrix algebra and give a new characterization of
inverse-closedness by means of finite sections. 
This is a necessary step in the qualitative analysis of the
convergence properties of the finite section method on weighted $\ell
^p$-spaces, but should be of independent interest in the study of
Banach algebras.

Let $\cA ^{FS} $ be the set of all finite sections of matrices in
$\cA $, formally 
\begin{equation}
  \label{eq:24}
\cA ^{FS} = \{ \vec{A} = (A_n)_{n\in \bN } : A_n = P_n A P_n \, \,
\mathrm{ for \, some } \, A \in \cA \} \, .
 \end{equation}
Although $\ac^{FS}$ is not
an algebra anymore, we may define a notion  of
inverse-closedness.  

\begin{definition}
We say that $\ac^{FS}$ is inverse-closed
if for every sequence $\{A_n\}_{n\in\bn}$ of invertible finite sections such that $\|A_n\|_\ac\le C$ and $\|A_n^{-1}\|_{op}\le C$, we have that 
$\|A_n^{-1}\|_\ac \le C'$, for some constants $C$ and $C'$ that do not depend on $n\in\bn$.
\end{definition}
The comparison of inverse-closedness  of $\cA $ and of $\cA ^{FS}$  indicates that the transition
from the infinite-dimensional setting $\ac$ to the finite-dimensional
case $\ac^{FS}$ is done by replacing the hidden word ``bounded" by
``bounded uniformly in dimension". 

The definition of $\cA ^{FS} $ suggests as a next step  to  consider
sequences of arbitrary finite square matrices instead of finite 
sections. To define a norm that is related to the $\cA $-norm, we
must assume that the norm $\|\cdot\|_\ac$ can be applied to arbitrary finite
matrices by defining an appropriate embedding of $(\Cst^{2n+1})^d$ into 
$\ell ^2 (\zd )$.
For $j\in \zd$ and $C_n = \{-n, \dots , n \}^d\subseteq \zd $ we
define an extension of the $(2n+1)^d \times (2n+1)^d $-matrix $B$  to
in infinite matrix $B^J$
of $B$, say $B^J=(B^J)_{k,l\in\bz}$ such that 
\begin{equation}\label{ext}
(B^J)_{j+k,j+l}=\begin{cases} (B)_{kl} & \text{for  } k,l \in C_n \\ 
0 & \text{otherwise}.
\end{cases}
\end{equation}
Then we define the norm of $B$ by 
\begin{equation}
  \label{eq:23}
\|  B\|_\ac =\|B^J\|_\ac \, .
\end{equation}
In the big picture, this definition makes sense only when the norm
does not depend on the embedding cube $J=j+C_n$. This requires an
additional property of $\cA $. 

Let $T_l, l\in \zd,$ denote the  translation operator $T_lf(k)=f(k-l)$
acting on $f\in\el$. 
We  say that the norm  of $\ac$ is {\it translation-invariant} if 
\begin{equation}\label{translation}
  \|T_{-l}AT_l\|_\ac=\|A\|_\ac \qquad \forall A\in \cA , l\in \zd \, .
\end{equation}
Clearly, if the norm of $\ac$ is translation-invariant,  then $\|B^J\|_\ac$
does not depend on the cube  $J=j+C_n$,  and we can apply $\|\cdot\|_\ac$ to finite
matrices. From now on, let us assume that $\|\cdot\|_\ac$ is
translation-invariant.

Similarly to~\eqref{eq:24} and analogous to~\cite[Section~1.2.2]{HRS01}
we introduce the set $\cA^{F}$ by 
\begin{equation*}
  \label{eq:25}
\cA ^{F} = \{ \vec{B} = (B_n)_{n\in \bN } : B_n \,\, \text{is a } \,
(2n+1)^d \times (2n+1)^d\,\, \text{matrix and}\,\, \sup _{n \in \bN } \|B_n \|_{\cA}
< \infty  \} \, , 
\end{equation*}
and we endow $\cA ^F $ with the norm
\begin{equation}
  \label{eq:26}
  \|\vec{B} \|_{\cA ^F} = \sup _{n \in \bN } \|B_n \|_{\cA} .
\end{equation}
If the norm of $\cA $ is translation-invariant, then $\cA ^F $ is a
well-defined object. We note that $\cA ^F $ is a Banach algebra 
contained in $\cB := \bigoplus _{n=1}^\infty \cB ((\bC ^{2n+1})^d)$.

\begin{definition}
We say that $\ac^F$ is inverse-closed
if for every sequence $\{B_m\}_{m\in\bz}$ of finite invertible matrices such that $\|B_m\|_\ac\le C$ and $\|B_m^{-1}\|_{op}\le C$, we have that 
$\|B_m^{-1}\|_\ac \le C'$, for some constants $C$ and $C'$ that do not depend on $m\in\bz$.
\end{definition} 

In view of Definition~\ref{definverse} this amounts to saying that the algebra $\cA ^F $ is
inverse-closed in  $\bigoplus _{n=1}^\infty  \cB ((\bC ^{2n+1})^d)$.

The  inverse-closedness of $\cA $,  $\cA ^{FS}$, and
$\cA ^F$  depends on the original algebra $\cA $. For the
study of the relations between them we introduce some further natural
conditions.

\begin{itemize}
\item[(C1)] Weak solidity: For every $A\in\ac$ there is a constant $C$
     such that $\|A_n\|_\ac\le C\|A\|_\ac$ for all $n\in\bn$. 
\item[(C2)] Weak inverse-closedness: For every $A\in\ac$ that is invertible 
     on $\el$ the condition $\sup _{n\in \bN } \|A_n^{-1}\|_\ac\le C$
     implies that $A^{-1}\in\ac$. 
\end{itemize}
Our third condition concerns an infinite matrix $B^{block}$ that is
built from blocks of finite square matrices $\{B_m\}_{m\in\bN}$ by
stacking them ``along the diagonal". For this we choose  a sequence
$j_m \in \zd $ such that 
sequence  of cubes  $ J_m = j_m +C_m \subset \bz$ is disjoint. Now set
\begin{equation}\label{block}
B^{block}=\sum_{m\in\bN} B^{J_m}_m,
\end{equation}
where $B^{J_m}_m$ is the extension of $B_m$ given in (\ref{ext}).
\begin{itemize}
\item[(C3)] Block norm equivalence: There exist constants $C,C'>0$,
  such that for every  $\vec{B} \in \cA ^F$
  \begin{equation}
    \label{eq:27}
 C\|B^{block}\|_\ac\le \sup_{m\in\bN}\|B_m\|_\ac\le
C'\|B^{block}\|_\ac
  \end{equation}
\end{itemize}
\begin{remark}\label{rem}
We want to point out  that in all settings $\ac$, $\ac^{FS}$, and
$\ac^F$ it suffices to show the inverse-closedness property for
positive matrices. Therefore, if necessary, one could restrict
conditions (C1)--(C3) to such matrices. We note that the upper bound
in \eqref{eq:27}   follows
already from condition (C1). 
\end{remark}
These three conditions are sufficient to show that the concepts of
inverse-closedness in $\ac$, $\ac^{FS}$ and $\ac^{F}$ are equivalent.

\begin{tm}\label{equiv}
Let $\ac\subset \Cal B(\ell ^2 (\zd ))$ be a Banach algebra such that
the norm of  $\ac$ is translation-invariant and  $\ac$ satisfies conditions
(C1)--(C3).  Then the following are equivalent:
\begin{itemize}
\item[a)] $\ac$ is inverse-closed in $\cB (\ell ^2 (\zd ))$.
\item[b)] $\ac^{FS}$ is inverse-closed.
\item[c)] $\ac^{F}$ is inverse-closed in $\cB :=\bigoplus _{n=1}^\infty 
\cB((\Cst^{2n+1})^d)$.
\end{itemize}
\end{tm}

\begin{proof}

c) $\Rightarrow$ b)  This implication  is clear, because each sequence
of finite sections $A_n$ belongs to  $\cA ^F$ by (C1). 

b) $\Rightarrow$ a) Since a matrix  $A $ is invertible, \fif\ the matrices $A^*A$ and $A
A^*$ are invertible, we may assume without loss of generality that $A$
is positive and invertible on $\ell ^2 (\zd )$.  So
 assume that $A\in\ac$ is positive and invertible on $\el$. We want to
show that, if $\cA ^{FS}$ is  inverse-closed, then $A^{-1}\in\ac$. 
 
By condition (C1), we have that $\|A_n\|_\ac\le C$ for all $n\in\bn$
and some constant  $C$ independent of $n\in\bn$.  Moreover, 
\eqref{eq:fe1} implies that $\|A_n^{-1}\|_{op}\le \|A^{-1}\|_{op}$.
 Since by assumption $\ac^{FS}$ is inverse-closed,  we
 obtain  that $\|A_n^{-1}\|_\ac\le C'$ for all $n\in\bn$ and some
 $C'>0$. By  condition (C2) we  obtain  $A^{-1}\in\ac$.

a) $\Rightarrow$ c)
We argue by contradiction and  show that if (c) fails, then so does
(a). Assume that $\ac^{F}$ is not inverse-closed. This means  that
there is a sequence of finite invertible matrices $\{B_m\}_{m\in\bN }$
such that $\|B_m\|_\ac$ and $\|B_m^{-1}\|_{op}$ are uniformly bounded
in $m\in\bN $, but $\sup_{m\in\bN}\|B_m^{-1}\|_\ac=\infty$.

Consider  the
corresponding matrix $B^{block}$  as given in (\ref{block}). Then  its
inverse  $(B^{block})^{-1}$ is a block matrix that corresponds to the
sequence $\{B_m^{-1}\}_{m\in\bz}$. Therefore,
$\|(B^{block})^{-1}\|_{op}=\sup_{m\in\bN}\|B_m^{-1}\|_{op}<\infty$, so
$B^{block}$ is invertible on $\el$. 
 Condition (C3) implies  that  $\|B^{block}\|_\ac\le
 C'\sup_{m\in\bz}\|B_m\|_\ac<\infty$, so $B^{block}\in\ac$. The same
 condition  guarantees that $\|(B^{block})^{-1}\|_\ac\ge
 C\sup_{m\in\bz}\|B_m^{-1}\|_\ac=\infty$. Thus 
 $(B^{block})^{-1}\notin\ac$ and $\ac$ cannot be 
 inverse-closed in $\cB (\el )$. 
\end{proof}

 Next we  apply  Theorem~\ref{equiv} to the matrix  algebras $\ac_v$,
 $\ac_v^1$, and $C_v$  introduced in Section~2.   As we
mentioned in Lemma~\ref{bound}, all these algebras are contained
in $\Cal B(\el)$. Moreover,  the norms associated to these algebras
are translation invariant. Indeed, to check that (\ref{translation}) holds, we
denote the standard basis of $\el$ by $\{e_k\}_{k\in\zd }$ and observe
that $\lan T_{-j}AT_j e_l,e_k\ran=\lan AT_j e_l,T_je_k\ran=\lan A
e_{l+j},e_{k+j}\ran$ for all $j\in \zd $. Therefore,
$(T_{-j}AT_j)_{kl}=(A)_{k+j,l+j}$ and 
since the norms of $\ac_v$,
 $\ac_v^1$ and $C_v$ use  only the difference of $k$ and
$l$, they are translation invariant.  

By Lemma~\ref{bound}(c)  each of these algebras is solid, so condition (C1)
holds for all of them. Condition (C3) is more problematic. Since the
norm of a matrix in  $\cA _v$ and $\cA ^1_v$  is defined in terms of its
rows and columns, it follows that    $\|B^{block}\|_\ac=
\sup_{m\in\bz}\|B_m\|_\ac$. So  property
(C3) holds for for  $\ac_v$  and  $\ac_v^1$.  Condition (C3) fails,
however, for $C_v$.

It remains to consider the weak inverse-closedness (C2). 

\begin{prop} 
Condition (C2) holds for each of the algebras   $\ac_v$ and   $\ac_v^1$.

\end{prop}
\begin{proof}
 Assume that $A\in\ac$ is invertible on $\el$ and that 
$ \sup _{n\in
   \bN } \|A_n^{-1}\|_\ac = C < \infty $.  
 Recall that  $\widetilde {A_n}= P_n AP_n + \lambda _+ (I-P_n)$ is the
extension of $A_n$  defined in \eqref{eq:1}.  Clearly, our assumption that $\|A_n^{-1}\|_\ac$ is uniformly
bounded, implies immediately that $\|\widetilde{A_n}^{-1}\|_\ac$ is
uniformly bounded as well. 

Since both $A$ and $A^{-1}$ are bounded on $\el$,
$\widetilde{A_n}^{-1}$ converges strongly to $A^{-1}$ in $\el$, as we
have seen in the proof of Theorem~\ref{cstar}, Step 3.  Therefore, $\lan
\widetilde{A_n}^{-1} e_l, e_k\ran$ 
converges to $\lan A^{-1} e_l, e_k\ran$ for all vectors of the
standard basis of $\el$.

\textbf{Case I: $\cA = \cA _v$.}  
Since 
$\|\widetilde{A_n}^{-1}\|_\ac$ is uniformly bounded, we have that 
\[
\|\widetilde{A_n}^{-1}\|_\ac=\sup_{k,l\in\bz}|\lan
\widetilde{A_n}^{-1} e_l, e_k\ran|v(k-l)\le C .
\]
Thus, we obtain that, for every $k,l\in\bz$, 
\[
|\lan A^{-1} e_l, e_k\ran|v(k-l)=\lim_{n\to\infty}|\lan \widetilde{A_n}^{-1} e_l, e_k\ran|v(k-l)\le C,
\]
and so $A^{-1}\in\ac$ with $\|A^{-1}\|_\ac\le C= \sup _{n\in \bN }
\|\widetilde{A_n}\inv \|_{\cA _v}$.

\textbf{Case II: $\cA = \cA _v ^1$.} We  use Fatou's Lemma. We have that 
\[
\|\widetilde{A_n}^{-1}\|_\ac=\max\bigg\{\sup_{k\in\bz}\sum_{l\in\bz}|\lan \widetilde{A_n}^{-1} e_l, e_k\ran|v(k-l),\sup_{l\in\bz}\sum_{k\in\bz}|\lan \widetilde{A_n}^{-1} e_l, e_k\ran|v(k-l)\bigg\}\le C.
\]
Therefore, we obtain, for every $k\in\zd$,  
\[
\sum_{l\in\bz}|\lan A^{-1} e_l, e_k\ran|v(k-l)\le \liminf_{n\to\infty}\sum_{l\in\bz}|\lan \widetilde{A_n}^{-1} e_l, e_k\ran|v(k-l)\le C
\]
and for every $l\in\bz$ 
\[
\sum_{k\in\bz}|\lan A^{-1} e_l, e_k\ran|v(k-l)\le \liminf_{n\to\infty}\sum_{k\in\bz}|\lan \widetilde{A_n}^{-1} e_l, e_k\ran|v(k-l)\le C.
\]
 Taking the  supremum over $l\in \zd $ (or $k$ respectively),  we conclude  
that $A^{-1}\in\ac^1_v$ and  $\|A^{-1}\|_{\ac _v^1}\le C= \sup _{n\in\bN }
 \|\widetilde{A_n}\inv \|_{\cA _v^1}$. 
\end{proof}

Since  we have verified  that all assumptions of Theorem \ref{equiv} are
 satisfied for $\ac_v$ and $\ac_v^1$, we obtain the following result.  
\begin{tm}\label{equiv2}
If $\ac$ is either $\ac_v$ or $\ac_v^1$, then $\cA $ is
inverse-closed, \fif\ $\cA ^{FS}$ is inverse-closed \fif\ $\cA ^{F}$ is
inverse-closed. 
\end{tm}
\begin{remark}
 Theorem~\ref{equiv2}  holds for $\ac_v^1$ even if $v\equiv1$. In this
case, $\ac_v^1$ is  the Schur class, i.e., the class of matrices that 
satisfy the Schur test~\cite{HLP52} or, equivalently, the class of matrices that
are bounded simultaneously  on all $\ell ^p, 1\leq p \leq
\infty $.  It seems to be  an open problem if
this algebra is inverse-closed.  Theorem~\ref{equiv2} reduces  this
problem to an equivalent (and equally difficult) question about
finite-dimensional matrices. However, if $v$ satisfies a mild growth
condition, see Theorem~\ref{inverseclosed}(b), then $\cA ^1_v$ is 
inverse-closed.  
\end{remark}

We need a few more facts before addressing questions
of convergence of the finite section method.

\begin{cor}\label{finsec}
Let $v$ be an admissible weight satisfying the
GRS-condition~\eqref{GRS}. 
 
(a) If $A\in \cA _v^1$ is positive  and invertible on $\el $, then 
$\sup_{n\in\bn}\|A_n^{-1}\|_{\cA _v^1} <\infty$.

(b) Assume in addition that   $v\inv \in \ell ^1(\zd )$ and  $v\inv \ast v\inv
\leq C v\inv $. If $A$ is positive and  invertible on $\el$, then
$\sup_{n\in\bn}\|A_n^{-1}\|_{\cA _v} <\infty$. 
\end{cor}
\begin{proof}
Our assumptions on $v$ imply that $\ac\in\{\ac_v,\ac_v^1\}$ is inverse-closed. By 
Theorem~\ref{equiv2}, $\ac^{FS}$ is inverse-closed as
well. Therefore, to achieve that   $\sup_{n\in\bn}\|A_n^{-1}\|_\ac<\infty$
it is enough to assure that $\|A_n\|_\ac$ and $\|A_n^{-1}\|_{op}$ are
bounded uniformly in $n\in\bn$. This, however, follows from
$\|A_n\|_\ac\le\|A\|_\ac$ (solidity) and from Step 1 in the proof of
Theorem~\ref{cstar}, 
where we showed that $\|A_n^{-1}\|_{op}\le\|A^{-1}\|_{op}$. 
\end{proof}

As the block norm equivalence (C3) fails for the algebra $\cC _v$, we
do not know whether  $(\cC _v)^{FS}$ is inverse-closed. However, since
$\cC _v \subseteq \cA _v ^1$, we have the following result. 

\begin{cor}
  \label{finsec2}
Let $v$ be an admissible weight satisfying the
GRS-condition~\eqref{GRS} and $v(k) \geq C (1+ |k|)^\epsilon $ for
some $\epsilon >0$. If  
 $A\in \cC _v $ is positive  and invertible on $\el $, then 
$\sup_{n\in\bn}\|A_n^{-1}\|_{\cA _v^1} <\infty$. 
\end{cor}

\section{Convergence in $\lpm $}
\label{convergence}

After the analysis of finite sections in matrix algebras,  we are now  in a
position to show that the
finite section method converges in weighted $\ell ^p$-spaces, whenever
the matrix is in one of the algebras 
 $\cA _s$, $\cA _v $,  $\avo $, and $\cC _v$. 



 



\begin{tm} \label{banachconv}
Let  $\cA $ be  one of the inverse-closed algebras $\cA _s, \cA _v , \cC _v
$, or $\cA _v^1$, where the  weight satisfies the conditions  stated
in Theorem~\ref{inverseclosed}  for each case. 
Assume that $A\in \cA $ is positive and invertible on $\el $ 
and acts boundedly  on $\lpm $. 

If $b\in \lpm$ and  $p< \infty $, then the finite section method
converges in the norm of $\lpm $.   

If $b\in \lpm$ and   $p = \infty $, then the finite section method
converges in the weak$^*$-topology. In particular, $x_n $ goes to $x$
entrywise.  
\end{tm}

\begin{proof}
We expand the model proof of Theorem~\ref{cstar} and insert the results
about Banach algebras obtained in Sections~3 and~4.  Recall that
$\widetilde {A_n}= P_n AP_n + \lambda _+ (I-P_n)$ is the 
extension of $A_n$  defined in~\eqref{eq:1}. Throughout the proof
$C$ denotes a constant that may change from step to step.

  \textbf{Step 1} in the proof of Theorem~\ref{cstar} remains
  unchanged and yields that  
$$
\sigma (\ant ) \subseteq [\lambda
  _-, \lambda _+]
$$
independent of $n$ and that 
\begin{equation}
  \label{eq:fe1a}
  \sup _{n\in \bN } \|\ant \inv \|_{op} \leq \lambda _-\inv = \|A\inv
  \|_{op} \, ,
\end{equation}
(where $\| \cdot \|_{op}$ is the operator norm on $\el $). 
Since $A $ is positive and invertible on $\el $ and $\cA $ is
inverse-closed in $\cB (\el )$ by our hypotheses on the weights,
Theorem~\ref{inverseclosed} guarantees that $A\inv \in \cA $ as well. By 
Corollary~\ref{closedgraph}(b), the inverse $A\inv $ is then bounded 
on $\lpm$. Furthermore, for $\cA \in \{ \cA_s, \cA_v, \cA_v^1\}$ by 
Corollary~\ref{finsec} we know that
$$
\sup _{n\in \bN } \|\ant \inv \|_{\lpm \to \lpm } \leq C \sup _{n\in
  \bN } \|\ant \inv \|_{\cA  } = C <\infty \, .
$$
For $\cA = \cC_v$,  Corollary~\ref{finsec2} implies that
$$
\sup _{n\in \bN } \|\ant \inv \|_{\lpm \to \lpm } \leq C \sup _{n\in
  \bN } \|\ant \inv \|_{\cA_v^1  } = C <\infty \, .
$$

\textbf{Step 2.}  For  $p< \infty $,  $\ant $ converges to $A$ in
the strong operator 
topology on $\lpm $. This follows from the inequality
\[
\|A-P_nAP_n \|_{\lpm \to \lpm } \leq \|(I-P_n)A \|_{\lpm \to\lpm}
+\|P_nA(I-P_n)\|_{\lpm \to \lpm },
\]
and the fact that $P_nf \to f \in \lpm$ is equivalent to the density
of the finite sequences in $\lpm$. 

\textbf{Step 3.} (Lemma of Kantorovich). We know that $\ant \inv = A_n
\inv + \lambda _+ \inv 
(I-P_n )$.  Since
\begin{eqnarray}
\|\ant \inv b - A\inv b \|_{\lpm}  &=&\| \ant \inv (A-\ant ) A\inv b \|_{\lpm }
\notag  \\
&\leq & \sup _n \|\ant \inv \|_{\lpm 
\to \lpm} \, \|  (A-\ant ) A\inv b \|_{\lpm }  \\
&\leq & \sup _n \|\ant \inv \|_{\cA} \, \|  (A-\ant ) A\inv b \|_{\lpm }
\label{kant1}
\end{eqnarray}
the strong convergence $\ant \rightharpoonup A $ on $\lpm $  implies that $\ant
\inv $ converges strongly to $A\inv $ on $\lpm $ for $1\le p < \infty$.

\textbf{Step 4.} Recall, that $A_n x_n = b_n $ and $Ax=b$. Then
\begin{eqnarray}
  \|x-x_n \|_{\lpm } &=& \|A\inv b - A_n \inv b_n \|_{\lpm } = \| A\inv b - A_n
  \inv P_n  b \|_{\lpm } \notag  \\
&\leq &\| (A\inv - \ant \inv ) b\|_{\lpm } + \| \ant \inv (b - P_n b)
\|_{\lpm } = \tI + \tII \, .
\end{eqnarray}
For $1 \le p < \infty$ the first term goes to zero by Step 3, and the 
second term is estimated by 
\begin{equation}
\tII \leq \sup _n \|\ant \inv \|_{\lpm \to \lpm } \, \|b-P_n b\|_{\lpm }
\leq C \inv \|b-P_n b\|_{\lpm },
\label{kant7}
\end{equation}
and also goes to zero. 

For $p=\infty$ we prove weak$^*$-convergence.  Assume
$b\in\ell^\infty_m=\Big(\ell^1_{1/m}\Big)^*$ and $y\in
\ell^1_{1/m}$. Then,  
\[
\lan x-x_n,y\ran=\lan A^{-1}b-\widetilde A_n^{-1}P_nb,y\ran=
\lan b-P_nb,A^{-1}y\ran+\lan P_nb,(A^{-1}-\widetilde A_n^{-1})y\ran.
\]
the first term tends to zero, because finite sequences are weak$^*$-dense in $\ell^\infty_m$. The second term is majorized by $\|b\|_{\ell_m^\infty}\|(A^{-1}-\widetilde A_n^{-1})y\|_{\ell^1_{1/m}}$ and converges to zero by Step 3, \eqref{kant1}.
\end{proof}

\section{Quantitative Estimates}
\label{quantitative}

In Theorem~\ref{banachconv} we have investigated the convergence of the finite
section method in the norm of $\lpm $ provided that the  input vector $b$
is in $ \lpm$.  For the quantitative analysis, we assume that the
input is in $\lpm$ and we study the convergence in a weaker norm. 

 We first  work with the algebra $\cA = \cA _v $ defined by off-diagonal
decay  of the matrices and a subconvolutive  weight  satisfying the
GRS-condition. Recall that $C_n = \{-n, \dots n \}^d$ is the cube of
integer vectors, so that $\sum _{k\not \in C_n} \dots $ becomes a tail
estimate.  



\begin{tm}
  \label{quant}
Assume that  $A\in \cA _v$ is invertible and $b\in \ell ^\infty _v (\zd )$. Set 
$$
\phi (n)  = \big( \sum _{k\not \in C_n} v(k)^{-2}\big) ^{1/2}\, .$$
 Then the finite section method
converges in the $\lz $-norm with the asymptotic estimate for the
error
\begin{equation}
  \label{eq:4}
  \|x-x_n \|_2 \leq C \phi (n) \, .
\end{equation}
\end{tm}

\begin{proof}
  We write 
  \begin{eqnarray*}
    x-x_n &=& A\inv b - A_n \inv b_n \\
&=& A\inv (b - b_n ) +  A\inv (A_n -A) A_n \inv P_n b = \tI + \tII \,.
  \end{eqnarray*}
Estimate of $\tI$: 
\begin{eqnarray*}
  \|\tI\|_2 \leq \|A\inv \|_{op} \,  \|b-b_n \|_2 &=& \|A\inv \|_{op} \,  \big(
  \sum _{|k|>n} |b_k|^2\big)^{1/2} \\
&\leq & \|A\inv \|_{op} \,  \, \|b\|_{\ell ^\infty _v} \big( \sum
_{k\not \in C_n}
v(k)^{-2}\big)^{1/2} \\ 
&=& \|A\inv \|_{op} \,  \, \|b\|_{\ell ^\infty _v} \,  \phi (n) \, .
\end{eqnarray*}
Estimate of $\tII$:
Set $z_n = AP_n A_n \inv P_n b= AP_n \ant \inv P_n b$, then $\tII =
A\inv (P_n - I)z_n$. Using
Lemma~\ref{bound}, Corollary~\ref{finsec} and the obvious fact that $\|P_n b\| _{\lpm } \leq
\|b\|_{\lpm }$ (true for every solid sequence space) we obtain 
$$
\|z_n \|_{\ell ^\infty _v} \leq \|A\|_{\cA _v} \|\ant \inv \|_{\cA _v}
\|b\|_{\ell ^\infty _v},$$
or the pointwise estimate
$$
|(z_n )_k| \leq C v(k) \inv 
$$
which is independent of $n$. 
So 
\begin{eqnarray*}
  \|\tII\|_2 &=& \|A\inv (P_n - I)z_n \|_2 \\
 &\leq & \| A\inv \|_{op} \,  \|(P_n -I) z_n\|_2 \\
&=& \| A\inv \|_{op} \,  \big( \sum _{k \not \in C_n} |(z_n )_k|^2\big)^{1/2} \\
&\leq & \| A\inv \|_{op} \,  C \big( \sum _{k \not \in C_n} v(k)^{-2}\big)^{1/2}
\\
&=& \| A\inv \|_{op} \,  C \phi (n) \, .
\end{eqnarray*}
\end{proof}
\begin{remark}
  If $v(x) = (1+|x|)^s$, then $\phi (n) = \Big(\sum _{k\not \in C_n}
  (1+|k||)^{-2s}\Big)^{1/2}  \sim n^{-s+d/2}$, and we
recover the result of ~\cite{Str00}. 
\end{remark}

Theorem~\ref{quant} can be generalized to other matrix algebras and
sequence spaces. For this we note that the input $b$ is in a ``small''
space, but that we measure the rate of convergence in a ``large''
space. The other item is that $A$ and $\widetilde A_n$ have to be
invertible on both the small and the large space with uniform bounds.

The rate of convergence will follow from the following tail estimates for the embedding of sequence spaces. 

\begin{lemma}\label{tail}
Assume that $\ell_m^p\subseteq\ell^q_w$ for $1\le p,q\le \infty$ and $m,w$ are moderate weights. Set $r^{-1}=\max\{q^{-1}-p^{-1},0\}$ and 
\begin{equation}
\label{rz1}
\varphi(n)=\Bigg( \sum_{k\notin C_n} \frac{w(k)^r}{m(k)^{r}}\Bigg)^{\frac1r}.
\end{equation}
Then, $\|b-P_nb\|_{\ell^q_w}\le\varphi(n)\|b\|_{\ell^p_m}$.
\end{lemma}

\begin{proof}
We write $b-P_nb=(1-\chi_{C_n})b$, where $\chi_{C_n}$ is the
characteristic function of $C_n$. Then,
$\|b-P_nb\|_{\ell^q_w}=\|bm(1-\chi_{C_n})\frac wm\|_{\ell^q}$. If
$p\le q$, then $\|c\|_{\ell^q}\le \|c\|_{\ell^p}$ and so  
\[
\|b-P_nb\|_{\ell^q_w}\le\Big\|bm(1-\chi_{C_n})\frac wm\Big\|_{\ell^p}
\le \Big\|(1-\chi_{C_n})\frac wm\Big\|_\infty\|bm\|_{\ell^p}
=\varphi(n)\|b\|_{\ell^p_m},
\]
(since for $r=\infty$ formula \eqref{rz1} has to be interpreted with  the
supremum norm). If $p>q$, then $r=(q^{-1}-p^{-1})^{-1}>1$. Thus, we
use H\"older's inequality $\|bc\|_q\le\|b\|_p\|c\|_r$ and obtain 
\[
\|b-P_nb\|_{\ell^q_w}\le\Big\|(1-\chi_{C_n})\frac wm\Big\|_r\|bm\|_{\ell^p}
=\varphi(n)\,\|b\|_{\ell^p_m},
\]
\end{proof}

\begin{tm}\label{errorquant}
Let $\cA$ be one of the inverse-closed algebras $\as,\av,\avo$ or
$\cC_v$. Assume that $\ell^p_m\subseteq\ell^q_w$ and that $\cA$ acts
boundedly on both $\ell^p_m$ and $\ell^q_w$. If
$A\in\cA$ is invertible on $\ell^2$ and $b\in\ell^p_m$ (the
``smaller'' space), then the finite section method converges in
$\ell^q_w$ (the ``larger'' space) with the error estimate 
\begin{equation}\label{errorest}
\|x-x_n\|_{\ell^q_w}\le C\|b\|_{\ell^p_m}\,\varphi(n),
\end{equation}
where $C=\|A^{-1}\|_{\ell^q_w}(1+\|A\|_{\ell^p_m}\,\|\widetilde
A_n^{-1}\|_{\cA})$ and $\varphi(n)$ is as in~\eqref{rz1}.
\end{tm}

\begin{proof}
As in the proof of Theorem \ref{quant} we estimate the error by 
\begin{equation}\label{xy}
\|x-x_n\|_{\ell^q_w}\le \|A^{-1}(b-b_n)\|_{\ell^q_w}+\|A^{-1}(P_n-I)z_n\|_{\ell^q_w},
\end{equation}
where $z_n=AP_n\widetilde A_n^{-1}P_nb$.

Since $A\in \cA$ is invertible on $\ell^2$, by inverse-closedness
$A^{-1}\in \cA$ and consequently $A^{-1}$ is also bounded on
$\ell^q_w$. We note that, by Corollaries~\ref{finsec} and
\ref{finsec2} we also have $\sup_{n\in\bn}\|\widetilde
A_n^{-1}\|_{\cA}<\infty$.  

For the first term in \eqref{xy} we obtain, with Lemma~\ref{tail}, that
\[
\|A^{-1}(b-b_n)\|_{\ell^q_w}\le \|A^{-1}\|_{\ell^q_w}\|b-b_n\|_{\ell^q_w}\le \|A^{-1}\|_{\ell^q_w}\|b\|_{\ell^p_m}\varphi(n).
\]

The second term is estimated by
$$
\|A^{-1}(P_n-I)z_n\|_{\ell^q_w}\le\|A^{-1}\|_{\ell^q_w \to \ell ^q_w}\,\|(P_n-I)z_n\|_{\ell^q_w}\le
\|A^{-1}\|_{\ell^q_w \to \ell ^q_w}\, \|z_n\|_{\ell^p_m}\, \varphi(n). 
$$

Finally, 
\[
\|z_n\|_{\ell^p_m}=\|AP_n\widetilde
A_n^{-1}P_nb\|_{\ell^p_m}\le\|A\|_{\ell^p_m \to \lpm }\|\widetilde
A_n^{-1}\|_{\lpm \to \lpm}\|b\|_{\ell^p_m}, 
\]
and this expression is uniformly bounded by Corollary~
\ref{finsec} and \ref{finsec2} and~Corollary~\ref{closedgraph}. Thus,
we are done. 
\end{proof}
If $\ell^\infty_v\subseteq\ell^2$, then we recover the simpler statement of Theorem~\ref{quant}.

\section{Non-Symmetric Finite Section Method for Non-Symmetric Matrices}
\label{s:nonsymm}

In the previous section we derived quantitative estimates for
the convergence of the finite section method under the assumption
that the matrices are positive definite.  This assumption
is crucial. For non-hermitian matrices it is already  a difficult
problem to derive merely qualitative statements about the convergence of
the finite section method~\cite{GF74,GGK03,RRS04,BS98}. Indeed, even
for very simple non-hermitian matrices the finite section method may fail. 

As an example, let us consider the Laurent operator given by the
biinfinite Toeplitz matrix 
$$
A = 
\begin{bmatrix}
\ddots &     &     &          &     &    &        \\
       & 0   & 0   & 0        & 0 & 0    &        \\
       & 1   & 0   & 0        & 0 & 0    &        \\
       & c   & 1   & \fbox{0} & 0 & 0    &        \\
       & c^2 & c   & 1        & 0 & 0    &        \\
       & c^3 & c^2 & c        & 1 & 0    &        \\
       &     &     &          &   &      & \ddots 
\end{bmatrix}.
$$
Here, we assume $|c|<1$ and, as usual in the finite section method literature,
the box indicates the entry in the zero-zero position.
An easy calculation shows that $A$ is invertible on $\ell^2(\Zst)$
and its inverse is the Laurent operator with biinfinite Toeplitz matrix
$$
A^{-1} = 
\begin{bmatrix}
\ddots &    &           &    &         \\
       & -c & 1         &  0 &         \\
       &  0 & \fbox{-c} &  1 &         \\
       &  0 &  0        & -c &         \\
       &    &           &    & \ddots 
\end{bmatrix}.
$$
Let us choose $b = e_0$, i.e, the right hand side is given by the
zero-th unit vector. The solution to $Ax=b$ is obviously the zero-th
column of $A^{-1}$, $x=e_{-1}+c e_0$. The finite section method as described
in~\eqref{Axb1}--\eqref{Axb2} applied to this system  fails
completely, because  none of the matrices $A_n$ is invertible. In theory the solution could be 
computed by solving the normal equations $\AT Ax = \AT b$. Thus, one might 
want to apply the finite section method to the positive definite 
system $A^{\ast}Ax = A^{\ast}b$ and invoke the results from the previous
sections, since $A$ (and thus $\AT$) belongs to $\cA_v,\cA_v^1$, or $\cC_v$ 
with weight $v(k)=e^{|k|^\alpha}, 0 < \alpha < 1$. However, the computation of 
$P_n \AT A P_n$ involves the infinite 
matrices $A,\AT$, which makes this approach not feasible for numerical
purposes.

 It is easy to see how to alter the
finite section method to make it work for this particular
example.  Our goal is  more ambitious, and  we want to derive a version
of the finite section method that  works for large classes of
(algebras of) non-hermitian matrices, and not just   for some individual
cases. We will  derive  conditions for the convergence of the 
finite section method for non-hermitian matrices in some matrix
algebras. 
For this, we consider a slightly generalized version of the finite
section method. 

\bigskip
Consider the system $Ax=b$ where $A$ is an invertible, but not necessarily 
hermitian matrix. We set
\begin{equation}
\Arn = P_r A P_n, \qquad \text{and} \qquad b_{r,n} = \Arn^{\ast} b,
\label{eq2}
\end{equation}
and try to solve the  system
\begin{equation}
\Arn^{\ast} \Arn \xrn = b_{r,n}
\label{eq3}
\end{equation}
for properly chosen $r$ and $n$.
Observe that $A_{r,n}$ is a $(2r+1)^d \times (2n+1)^d$ matrix, and so
$A_{r,n}^* A_{r,n}$ is a $(2n+1)^d \times (2n+1)^d$-matrix. 
 In general we will need $r>n$, therefore
we refer to~\eqref{eq2}--\eqref{eq3} as {\em non-symmetric finite section
method}.

Let us denote $B_n:=P_n A^{\ast} A P_n$,
$D_{r,n}:= P_n A^{\ast} P_r A P_n =A_{r,n}^{\ast} A_{r,n}$. 
Analogously to~\eqref{eq:1} we define the extensions
\begin{equation}
\widetilde{B_n} = B_n + \lambda_{+} (I- P_n),\qquad
\widetilde{D_{r,n}} = D_{r,n} + \lambda_{+} (I- P_n),
\label{ext1}
\end{equation}
where $\sigma(A^{\ast}A) \subseteq [\lambda_{-},\lambda_{+}]$.

Clearly, $B_n, n\in \bN ,$ is the sequence of finite sections of
$A^*A$ and $D_{r,n}$ is an approximation of $B_n$. We  study this
approximation for matrices in $\cA _v$.

\begin{lemma}\label{convnonsym}
  Assume that $A\in \cA _v$. Then there exists a sequence $R(n) \in
  \bN $, such that for every $r(n) \geq R(n)$ 
  \begin{equation}
    \label{eq:ll1}
    \lim _{n\to \infty } \|B_n - D_{r(n),n} \|_{\av } = 0 \, .
  \end{equation}
If $v(k) = (1+|k|)^s$ and $\av = \cA _s$, then we may choose  $R(n)= n^\alpha$
for $\alpha > \frac{2s}{2s-d}$ and obtain the rate  
$$
\| B_n - D_{n^\alpha ,n} \|_{\cA _s} \leq C \|A\|_{\cA _s} \, n^{\alpha (d-2s) +
  2s} \, , 
$$
\end{lemma}

\begin{proof}
We define $\TEr = \TBn - \TCr= B_n - D_{r,n }$.
Clearly, $\TEr$ is hermitian and in $\cA _v$, and 
$(\widetilde{E_{r,n}})_{kl} = 0$ for $k,l \not \in C_n $. 
If  $k,l\in C_n$, then 
$$(E_{r,n})_{kl}=  (B_n)_{kl}-(D_{r,n})_{kl} = \sum _{j\in \zd }
(A^*)_{kj} a_{jl} - \sum_{j \in C_r} (A^*)_{kj} a_{jl} = \sum _{j \not
  \in C_r}\overline{a_{jk}} a_{jl} \, ,$$
and we obtain the estimate
$$
|(\TEr )_{kl} | \leq \ \sum _{j \not
  \in C_r}|a_{jk}|| a_{jl}| \, $$
for all entries. If $A\in \av $, then $|a_{jk}| \leq \|A\|_{\av }
v(j-k)\inv $. Consequently we estimate the norm of $\TEr$ by 
\begin{eqnarray*}
  \|\TEr \|_{\av } &=& \sup _{k,l \in C_n} |(\TEr )_{kl} | v(k-l) \\
  &\leq & \sup _{k,l \in C_n} \|A\|_{\av } ^2  \sum _{j\not \in C_r} 
v(j-k)^{-1}v(j-l)^{-1}v(k-l) \, .
\end{eqnarray*}
Since $v(j-l)\inv \leq v(j-k)\inv v(k-l)$, we continue with 
$$
  \|\TEr \|_{\av } \leq  \|A\|_{\av } ^2  \, \sup _{k,l \in C_n}  \sum
  _{j\not \in C_r}  v(j-k)^{-2} v(k-l)^2 \, . 
$$
Clearly, if $k,l \in C_n$ and $j\not \in C_r$, then $k-l \in C_{2n}$
and $j-k \not \in C_{r-n}$ and we arrive at the estimate
\begin{equation}
  \label{eq:ll2}
  \|\TEr \|_{\av } \leq  \|A\|_{\av } ^2  \, \sup _{k \in C_{2n}}
  v(k)^2 \,  \sum
  _{j\not \in C_{r-n}}  v(j)^{-2}  \, . 
\end{equation}
As a consequence, we obtain that $\lim _{r\to \infty } \|\TEr \|_{\av
} =0$ and \eqref{eq:ll1} is proved. 

If $A\in \cA _s$, i.e.,  $v(k) = (1+|k|)^s$, then $\sup _{k \in
  C_{2n}} v(k)^2 = \cO( n^{2s})$ and $ \sum
  _{j\not \in C_{r-n}}  v(j)^{-2} = \cO ( (r-n)^{d-2s})$. For $r(n) =
  n^\alpha $,  we
  obtain the explicit estimate
$$
\| B_n - D_{n^\alpha , n} \|_{\cA _s} \leq C \|A\|_{\cA _s} \, n^{\alpha (d-2s) +
  2s} \, , 
$$
which tends to $0$ for $\alpha > \frac{2s}{2s-d}$. 
\end{proof}

\begin{tm}
\label{th:nonsymm}
Let $\cA \in \{\as,\av\}$ where the weight $v$ satisfies the
conditions stated in Theorem~\ref{inverseclosed}(a).  Let $Ax=b$ be
given.  Assume that $b\in \lpm $ and that  $A \in \cA$ is invertible 
on $\ell^2(\Zst^d)$ and acts on $\ell^p_m$. 

Consider the finite sections
\begin{equation}
\Arn^{\ast} \Arn x_{r,n} = \Arn^{\ast}b.
\label{normalfin}
\end{equation}
Then, for every $n$ there exists an $R(n)$  (depending on $\lambda_{-}$ 
and $v$)  such that $x_{r(n),n}$ converges to $x$ in the norm of
$\ell^p_m$, for every choice  $r(n)\ge R(n)$.
\end{tm}

\begin{proof}
We split the error $x-x_{r,n}$ into three terms as follows:
\begin{align}
\|x - \xrn \|_{\ell^p_m} & = \|(\AT A)^{-1}\AT b - D_{r,n}\inv \Arn b \|_{\ell^p_m} 
\notag \\
& \le \|(\AT A)^{-1}\AT b - B_n^{-1}P_n \AT b\|_{\ell^p_m} +
 \|B_n^{-1}P_n \AT b - B_n^{-1} \Arn b\|_{\ell^p_m} +  \notag \\  
 & \,\, +    \|B_n^{-1}\Arn b -D_{r,n}\inv  \Arn b\|_{\ell^p_m} =
 \|\tI\|_{\lpm}+\|\tII\|_{\lpm}+\|\tIII\|_{\lpm} \, .
\label{est1}
\end{align}
We observe that the vector $B_n\inv P_n A^*b$ is exactly the result of
the finite section method applied to the normal equation $A^*Ax =
A^*b$. Since $A^*A \in \cA _v$ and $A^*b\in \lpm $, 
Theorem~\ref{banachconv} is applicable and implies that $\|\tI\|_{\lpm}  \to 0$ for
$p<\infty $ and $\tI \to 0$ weak$^*$ for $p=\infty $.

Since $\Arn = P_r A P_n$ and $B_n \inv P_n = \widetilde{B_n}\inv P_n$
we can  estimate the second term  by 
\begin{eqnarray}
  \|\tII\|_{\lpm } = \|\widetilde{B_n}\inv (P_n A^* b - P_n A^* P_r b)
  \|_{\lpm } \notag \\
\leq \sup _{n\in \bN} \|\widetilde{B_n}\inv  \|_{\lpm \to \lpm  } \,
\|A^*\|_{\lpm } \|b-P_r b\|_{\lpm } \, . \label{br}
\end{eqnarray}
As in the proof of Theorem~\ref{banachconv}, 
Corollary~\ref{closedgraph} and~\ref{finsec} imply  that 
$\sup _{n\in \bN} \|\widetilde{B_n}\inv  \|_{\lpm \to \lpm  } \leq  C\sup
_{n\in \bN} \|\widetilde{B_n}\inv  \|_{\av } \leq C' <\infty $. Since 
the finite sequences are dense in  $\lpm $ for $p<\infty $, \eqref{br}
yields $\|\tII\|_{\lpm } \to 0$, similarly $\tII \to 0$ weak$^*$ for
$p=\infty $. 

For the third term, we start with  the obvious  estimate
$$
\|\tIII\|_{\lpm } = \| B_n^{-1}\Arn b -\Dri \Arn b\|_{\ell^p_m}\leq
\|\widetilde{B_n}\inv - \widetilde{\Dri} \|_{\lpm \to \lpm } \,
\|A_{r,n}^* b\|_{\lpm } \, .
$$
Here $   \|A_{r,n}^* b\|_{\lpm }  = \|P_n A^* P_r b\|_{\lpm } \leq
\|A^*\|_{\lpm \to \lpm } \, \|b\|_{\lpm }$ is uniformly bounded
independent of $n$ and $r$.

For the operator norm we use an estimate for inverses in Banach
algebras, see e.g., ~\cite{conway90}, and obtain that 
\begin{eqnarray*}
  \|\widetilde{B_n}\inv - \widetilde{D_{r,n}}\inv  \|_{\lpm \to \lpm } &\leq &
C  \|\widetilde{B_n}\inv - \widetilde{D_{r,n}}\inv  \|_{\av} \\
&\leq & C \, \frac{\|\widetilde{B_n}\inv \|_{\av }^2 \, \|\widetilde{B_n} -
  \widetilde{D_{r,n}}\|_{\av }}{1-\|\widetilde{B_n}\inv \|_{\av }\, \|\widetilde{B_n} -
  \widetilde{D_{r,n}}\|_{\av } } \, .
\end{eqnarray*}
Once again by Corollary~\ref{finsec} we have $\sup
_{n\in \bN} \|\widetilde{B_n}\inv  \|_{\av } \leq C <\infty $, and by
Lemma~\ref{convnonsym} 
$    \lim _{r\to \infty } \|B_n - D_{r,n} \|_{\av } = 0 $. 

Consequently, for any  positive sequence $\epsilon _n \to 0$, we may choose $R(n)$,
such that $\|B_n - D_{r(n),n} \|_{\av } < \epsilon _n$ for $r(n) \geq R(n)$.

By combining the estimates for $\tI, \tII, \tIII$, we have thus proved that
$\|x- x_{r(n),n}\|_{\lpm } \to 0$ for every sequence $r(n) \geq R(n)$
and we are done. 
 \end{proof}

\begin{remark}
If $v(k) = (1+|k|)^s$ for $s>d$,  then $R(n)$ can be chosen to be
$n^\alpha $ for $\alpha >\frac{2s}{2s-d}$ by Lemma~\ref{convnonsym}.
\end{remark}
\begin{remark}
It is well-known that, from a numerical viewpoint, the  solution of
the  normal equations should be avoided whenever the 
condition number of the matrix is large.  As an alternative to the
normal equations one could  use matrix factorization
methods. Since  $D_{r,n}$ is invertible, the matrix $\Arn$ has full
rank $(2n+1)^d$, and  one  could  apply
a QR-factorization of  $\Arn$ or some other factorization and  compute
an approximate solution 
to $Ax=b$ in that way. This idea  raises a number of interesting questions:
For instance, assume we can factorize a matrix $A \in \cA$ into $A= QR$,
where $Q$ is unitary and $R$ is upper triangular, do the individual 
components $Q$ and $R$ also belong to $\cA$? How about other
matrix factorizations such as $LU$- or polar-decomposition?
We refer the reader to~\cite{Str06} for answers to these questions.
\end{remark}

\bigskip
We return to the example in the beginning of this section. Clearly,
$A$ belongs to  $\cA_v$ for every 
weight $v(k) = e^{|k|^{\alpha}}, 0 < \alpha < 1$. Since  the entries
of $A$ decay exponentially off the diagonal, it is not difficult to see
that it is sufficient to choose  $r(n) =sn$ for a sufficiently
large $s>1$, independently of $n$. In this particular example it would
even suffice to set $s=n+1$, but as pointed out, our goal was to derive
a finite section technique that is applicable to large classes of matrices,
not just to this particular one.

In light of Theorem~\ref{th:nonsymm} it is worthwhile to recall that a 
necessary and sufficient condition for the applicability of the finite section 
method~\eqref{Axb1}--\eqref{Axb2} to Laurent operators is that 
the winding number of the invertible Laurent operator is zero, 
cf.~\cite{GGK03,HRS01}. For the non-symmetric finite section method the 
winding number is not relevant, the key property is the off-diagonal
decay of the matrix.
Thus Theorem~\ref{th:nonsymm} considerably enlarges the range of 
applicability of finite section type methods even for the classical
and thoroughly analyzed cases of Laurent and Toeplitz operators.

\section{An example from digital communication} \label{s:example}

In this section we demonstrate the practical relevance of the
theoretical framework derived in this paper by analyzing a
problem arising in digital communication. We highlight the details
related to the finite section method and refer the reader
to~\cite{Pro00} for a more detailed description 
of the engineering aspects of the problem.

In a time-invariant digital communication system one is confronted with 
a linear system of equations $Ax=b$, where $x=\{x_l\}_{l\in\Zst}$ is a
sequence of information symbols to be transmitted and $b=\{b_k\}_{k\in\Zst}$
is the received, discrete signal. We can assume that
$x_k \in \{-1,1\}$, thus $x \in \ell^{\infty}$.
The entries of $A$ are of the form
\begin{equation}
a_{kl} = \phi(\cdot - k R) \ast h \ast \phi(\cdot - lT),
\label{commeq}
\end{equation}
where $h$ is the {\em channel impulse response}, $\phi$ is a bandlimited
function (the {\em transmission pulse}), $T$ is the transmission period, and
$R$ is the receive sampling period. We do not go into detail about the
particular choice of $\phi$, $T$, and $R$. The only facts we need are:
(i) for properly selected $T$ and $R$ we can choose $\phi$ to be a bandlimited
function in $L^1_v(\Rst)$, where $v$ must satisfy the Beurling-Domar condition,
i.e., 
\begin{equation}
\sum_{k=0}^{\infty} \frac{\log v (kx)}{k^2} < \infty, \qquad \text{for all
$x \in \Rst$;}
\label{BD}
\end{equation}
(ii) under certain conditions on $h$, the matrix $A$ has an inverse for $R=T$
and a left-inverse for $R < T$.

Furthermore, we note that $h$ is a causal function that decays exponentially 
in time. This implies that $A$ is non-hermitian and that $A \in \av$, 
the latter follows from well-known properties of Beurling convolution 
algebras~\cite{RS00} and the fact that a weight which satisfies~\eqref{BD}
also satisfies the GRS condition~\eqref{GRS}, cf.~\cite{Gro06a}.

There are two ways to approach the problem of recovering $x$ from $b$.
In the first case we try to recover the entries of $x$ ``on the fly'', 
i.e., we solve the truncated system $A_{m,n} x_n = b_n$. In this case we only
assume that $b \in \ell^{\infty}$ and the $\ell^p_m$-convergence estimates
of Theorem~\ref{th:nonsymm} apply.
In the second case we precompute the inverse of $A$ by solving
$Az=e_0$ where $e_0$ is the zeroth unit vector. Due to the specific
(block)-Toeplitz structure of $A$, the vector $z$ contains all
required information to fully determine the inverse of $A$, which
is then used to recover $x$.
In this case we apply the non-symmetric finite section method from
Section~\ref{s:nonsymm} to $Az = e_0$. Since $A \in \av$ and $A^{\ast}e_0 \in
\ell^1_v$, quantitative estimates as in Section~\ref{quantitative} apply and 
we can approximate the true solution $z$ with a rate of convergence
depending on $v$. Since in this application $v$ can be chosen to be 
$v(x) = e^{|x|^\alpha}$ with $\alpha <1$, the (non-symmetric) finite
section method achieves exponential rate of convergence.

\def\cprime{$'$} \def\cprime{$'$} \def\cprime{$'$} \def\cprime{$'$}
  \def\cprime{$'$}



\begin{thebibliography}{10}

\bibitem{Arv94}
W.~Arveson.
\newblock {$C\sp *$}-algebras and numerical linear algebra.
\newblock {\em J. Funct. Anal.}, 122(2):333--360, 1994.

\bibitem{Bas90}
A.~G. Baskakov.
\newblock Wiener's theorem and asymptotic estimates for elements of inverse
  matrices.
\newblock {\em Funktsional. Anal. i Prilozhen.}, 24(3):64--65, 1990.

\bibitem{Bas97}
A.~G. Baskakov.
\newblock Estimates for the elements of inverse matrices, and the spectral
  analysis of linear operators.
\newblock {\em Izv. Ross. Akad. Nauk Ser. Mat.}, 61(6):3--26, 1997.

\bibitem{BS83}
A.~B{\"o}ttcher and B.~Silbermann.
\newblock The finite section method for {T}oeplitz operators on the
  quarter-plane with piecewise continuous symbols.
\newblock {\em Math. Nachr.}, 110:279--291, 1983.

\bibitem{BS90}
A.~B{\"o}ttcher and B.~Silbermann.
\newblock {\em Analysis of {T}oeplitz operators}.
\newblock Springer-Verlag, Berlin, 1990.

\bibitem{BS98}
A.~B{\"o}ttcher and B.~Silbermann.
\newblock {\em Introduction to large truncated {T}oeplitz matrices}.
\newblock Universitext. Springer-Verlag, New York, 1999.

\bibitem{CS04}
O.~Christensen and T.~Strohmer.
\newblock The finite section method and problems in frame theory.
\newblock {\em Journal Approx.\ Theory}, 133(2):221--237, 2005.

\bibitem{conway90}
J.~B. Conway.
\newblock {\em A course in functional analysis}.
\newblock Springer-Verlag, New York, second edition, 1990.


\bibitem{DMS84}
S.~Demko, W.~F. Moss, and P.~W. Smith.
\newblock Decay rates for inverses of band matrices.
\newblock {\em Math. Comp.}, 43(168):491--499, 1984.

\bibitem{GRS64}
I.~Gelfand, D.~Raikov, and G.~Shilov.
\newblock {\em Commutative normed rings}.
\newblock Chelsea Publishing Co., New York, 1964.
\newblock Translated from the Russian.

\bibitem{GF74}
I.~Gohberg and I.~Fel'dman.
\newblock {\em Convolution equations and projection methods for their
  solution}.
\newblock American Mathematical Society, Providence, R.I., 1974.
\newblock Translated from the Russian by F. M. Goldware, Translations of
  Mathematical Monographs, Vol. 41.

\bibitem{GGK03}
I.~Gohberg, S.~Goldberg, and M.~Kaashoek.
\newblock {\em Basic classes of linear operators}.
\newblock Birkh\"auser Verlag, Basel, 2003.

\bibitem{GKW89}
I.~Gohberg, M.~A. Kaashoek, and H.~J. Woerdeman.
\newblock The band method for positive and strictly contractive extension
  problems: an alternative version and new applications.
\newblock {\em Integral Equations Operator Theory}, 12(3):343--382, 1989.

\bibitem{Gro06a}
K.~Gr\"ochenig.
\newblock Weight functions in time-frequency analysis, 2006.
\newblock preprint.

\bibitem{GL03}
K.~Gr{\"o}chenig and M.~Leinert.
\newblock Wiener's lemma for twisted convolution and {G}abor frames.
\newblock {\em J. Amer. Math. Soc.}, 17:1--18, 2004.

\bibitem{GL04}
K.~Gr{\"o}chenig and M.~Leinert.
\newblock Symmetry of matrix algebras and symbolic calculus for infinite
  matrices.
\newblock {\em Trans.\ Amer.\ Math.\ Soc.}, 358:2695--2711, 2006.

\bibitem{HRS01}
R.~Hagen, S.~Roch, and B.~Silbermann.
\newblock {\em {$C\sp *$}-algebras and numerical analysis}, volume 236 of {\em
  Monographs and Textbooks in Pure and Applied Mathematics}.
\newblock Marcel Dekker Inc., New York, 2001.

\bibitem{HLP52}
G.~H. Hardy, J.~E. Littlewood, and G.~P{\'o}lya.
\newblock {\em Inequalities}.
\newblock Cambridge, at the University Press, 1952.
\newblock 2d ed.

\bibitem{Jaf90}
S.~Jaffard.
\newblock Propri\'et\'es des matrices ``bien localis\'ees'' pr\`es de leur
  diagonale et quelques applications.
\newblock {\em Ann. Inst. H. Poincar\'e Anal. Non Lin\'eaire}, 7(5):461--476,
  1990.

\bibitem{RRS04}
V.~Rabinovich, S.~Roch, and B.~Silbermann.
\newblock {\em Limit operators and their applications in operator theory},
  volume 150 of {\em Operator Theory: Advances and Applications}.
\newblock Birkh\"auser Verlag, Basel, 2004.

\bibitem{Pro00}
J.~G. Proakis. 
\newblock {\em Digital Communications}.
\newblock McGraw-Hill, New York, 2000.

\bibitem{RRS01}
V.~S. Rabinovich, S.~Roch, and B.~Silbermann.
\newblock Algebras of approximation sequences: finite sections of
  band-dominated operators.
\newblock {\em Acta Appl. Math.}, 65(1-3):315--332, 2001.
\newblock Special issue dedicated to Antonio Avantaggiati on the occasion of
  his 70th birthday.

\bibitem{RS00}
H.~Reiter and J.~Stegeman.
\newblock {\em Classical harmonic analysis and locally compact groups},
  volume~22 of {\em London Mathematical Society Monographs. New Series}.
\newblock The Clarendon Press Oxford University Press, New York, second
  edition, 2000.

\bibitem{Sjo95}
J.~Sj{\"o}strand.
\newblock Wiener type algebras of pseudodifferential operators.
\newblock In {\em S\'eminaire sur les \'Equations aux D\'eriv\'ees Partielles,
  1994--1995}, pages Exp.\ No.\ IV, 21. \'Ecole Polytech., Palaiseau, 1995.

\bibitem{Str98a}
T.~Strohmer.
\newblock Rates of convergence for the approximation of dual shift-invariant
  systems in {$\ell^2({\mathbb Z})$}.
\newblock {\em J.\ Four.\ Anal.\ Appl.}, 5(6):599--615, 2000.

\bibitem{Str00}
T.~Strohmer.
\newblock Four short stories about {T}oeplitz matrix calculations.
\newblock {\em Linear Algebra Appl.}, 343/344:321--344, 2002.
\newblock Special issue on structured and infinite systems of linear equations.

\bibitem{Str06}
T.~Strohmer.
\newblock Matrix factorizations and {B}anach algebras, 2006.
\newblock manuscript.

\bibitem{Str05}
T.~Strohmer.
\newblock Pseudodifferential operators and {B}anach algebras in mobile
  communications.
\newblock {\em Applied and Computational Harmonic Analysis},
   20(2):237--249,2006.


\bibitem{Sun05}
Q.~Sun.
\newblock Wiener's lemma for infinite matrices with polynomial off-diagonal
  decay.
\newblock {\em C. R. Math. Acad. Sci. Paris}, 340(8):567--570, 2005.

\end{thebibliography}

\end{document}